\documentclass[12pt,oneside]{amsart}
\usepackage{amsthm,amsfonts,amssymb,amscd,euscript}

\usepackage[matrix,arrow]{xy}

\author{Florin Ambro} 
\address{RIMS, Kyoto University\\
Kyoto 606-8502, Japan.}
\email{ambro@kurims.kyoto-u.ac.jp}

\newcommand{\Q}{{\mathbb Q}}
\newcommand{\Z}{{\mathbb Z}}
\newcommand{\N}{{\mathbb N}}
\newcommand{\R}{{\mathbb R}}

\newcommand{\cB}{{\mathcal B}}
\newcommand{\cO}{{\mathcal O}}

\newcommand{\emb}{\operatorname{emb}}
\newcommand{\Exc}{\operatorname{Exc}}
\newcommand{\Int}{\operatorname{int}}
\newcommand{\relint}{\operatorname{relint}}
\newcommand{\lct}{\operatorname{lct}}
\newcommand{\Lct}{\operatorname{Lct}}
\newcommand{\Mld}{\operatorname{Mld}}
\newcommand{\mult}{\operatorname{mult}}
\newcommand{\Proj}{\operatorname{Proj}}

\theoremstyle{plain}
\newtheorem{thm}{Theorem}[section]

\newtheorem{lem}[thm]{Lemma}

\newtheorem{conj}[thm]{Conjecture}

\theoremstyle{definition}
\newtheorem{defn}[thm]{Definition}

\newtheorem{exmp}[thm]{Example}

\theoremstyle{remark}

\begin{document}

\bibliographystyle{amsalpha+}
\title[The minimal log discrepancy]
{The minimal log discrepancy}
\maketitle

\begin{abstract}
We survey the known and expected properties of the
minimal log discrepancy, the local invariant of a
log variety.  
\end{abstract}


\tableofcontents



\section*{Introduction}

\footnotetext[1]{The author is supported by a 
21st Century COE Kyoto Mathematics Fellowship,
and by the JSPS Grant-in-Aid No 17740011.
}

 The minimal log discrepancy is a fundamental invariant of
the singularities that appear in the birational classification 
of algebraic varieties. 
Introduced by Shokurov~\cite{Sh88} in connection to the 
termination of a sequence of flips, it appears in the local 
context of the classification of singularities, or the global 
context of Fujita's conjecture 
on adjoint linear systems. The minimal log discrepancy
measures stable vanishing orders of sections
of canonical graded rings, or the rate of growth of certain 
subspaces in the spaces of jets of a singularity. It has an 
arithmetic flavour, being related to the first 
minimum of Minkowski in the geometry of numbers. 

In this note we introduce the minimal log discrepancy,
present some basic open problems, and illustrate 
them with toric examples. We hope to reinforce the 
original connection of Reid~\cite{Reid79} between 
discrepancies of singularities on the one hand, 
and stable vanishing orders of sections of canonical
graded rings on the other hand.

The minimal log discrepancy is the local invariant 
of a {\em log variety}. We recall in \S 1 the construction 
of canonical models and discrepancies, and their logarithmic 
version. This is the natural motivation for log varieties 
with {\em log canonical singularities}, which locally  
are just open subsets of log canonical models.
We give the rigorous definition of log varieties and 
minimal log discrepancies in \S 2, and present explicit 
combinatorial formulas in the toric case. We present problems 
on minimal log discrepancies in \S 3, and discuss their 
toric case and some methods, old and new. In \S 4, we 
discuss the log canonical threshold and its applications 
to the problem of constructing {\em flat log structures}. 


\section{Background}



\subsection{Zariski decomposition}

The origin of Zariski decomposition is the problem
of computing the Hibert function
$\varphi(n)=\dim_{\mathbb C} R(X,D)_n$
of the graded ring 
$$
R(X,D)=\bigoplus_{n=0}^\infty H^0(X,nD),
$$
for a given divisor $D$ on a complex manifold $X$~\cite{Zariski62}.
We recall the solution to this problem in the 
case when $R(X,D)$ is finitely generated and $D$ is big.
We use Zariski's notation
$$
H^0(X,D)=\{a\in {\mathbb C}(X)^\times; (a)+D\ge 0\}\cup\{0\},
$$
which makes sense even if $D$ does not have integer 
coefficients.

First, if $D$ is an ample divisor, there exists a 
polynomial $P\in \Q[T]$ of degree $d=\dim(X)$, 
and a positive integer $n_0$, such that $\varphi(n)=P(n)$ 
for $n\ge n_0$. In general, finite generation means that 
there exists a positive integer $r$ such that the natural 
map $S^lH^0(X,rD)\to H^0(X,lrD)$ is surjective for 
every $l\ge 1$. By Hironaka's resolution of the base 
locus of a linear system, there exists a birational 
modification $f_r\colon X_r\to X$ such that if $F_r$ 
is the fixed divisor of $\vert f_r^*(rD)\vert$, the mobile
part $M_r=f_r^*(rD)-F_r$ defines a linear system without base
points. After taking the Stein factorization, $\vert M_r\vert$
defines a morphism with connected fibers $g_r\colon X_r\to Y_r$. 
The variety $Y_r$ has normal singularities and $M_r=g_r^*(A_r)$,
for a normally generated ample Cartier divisor $A_r$ on $Y$.
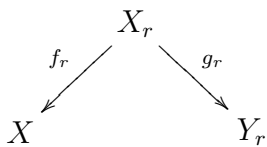
\[ \xymatrix{
& X_r \ar[dl]_{f_r} \ar[dr]^{g_r} &  \\
X  & & Y_r
}\]
We have $R(X,D)=R(X_r,f_r^*(D))$, the natural 
inclusion $R(X_r,\frac{1}{r}M_r)\subseteq R(X_r,f_r^*D)$ 
becomes an identity, and 
$R(X_r,\frac{1}{r}M_r)=R(Y,\frac{1}{r}A_r)$. In
particular, $Y_r=\Proj R(X,D)$ and 
$$
\varphi(n)=\dim_{\mathbb C}H^0(Y_r,\frac{n}{r}A_r).
$$ 
There are polynomials $P_0,\ldots,P_{r-1}\in \Q[T]$, 
of degree $d$, and a positive integer $n_0$, such that 
$\varphi(n)=P_{n-r\lfloor \frac{n}{r}\rfloor}(n)$ for
$n\ge n_0$. The growth rate of $\varphi$ is
$$
\lim_{n\to \infty}\frac{\varphi(n)}{n^d/d!}=\frac{(A_r^d)}{r^d}
\in \frac{1}{r^d}\N.
$$

The decomposition $f_r^*(D)=g_r^*(\frac{1}{r}A)+\frac{1}{r}F_r$
is the prototype of a Zariski decomposition. 
In general, we say that $D=P+F$ is a {\em Zariski decomposition} 
on a normal
variety $X$, if $F$ is an effective $\R$-divisor and $P$ is a
nef $\R$-divisor such that the induced inclusion 
$R(X,P)\subseteq R(X,D)$ is surjective. 
The Zariski decomposition is unique if $D$ is big, and 
in this case $F=\lim_{n\to \infty} \frac{1}{n}F_n$, 
where $F_n$ is the fixed part of $\vert nD\vert$.
There are several notions of Zariski decomposition 
at present, but they all coincide for big divisors
(see~\cite{Zariski62, Kawamata85, Fujita86, Cutkosky86, 
Nakayama04}).

Two useful examples of Zariski decomposition are
as follows. For a nef $\R$-divisor $P$,
the decomposition $P=P+0$ is a Zariski decomposition.
If $D=P+F$ is a Zariski decomposition on $X$ and 
$F'$ is an effective $\R$-divisor supported by the 
exceptional locus of a birational map $\mu\colon X'\to X$, 
then $\mu^*D+E=\mu^*P+(\mu^*F+F')$ is a Zariski decomposition.


\subsection{Canonical models, discrepancies}

Let $X$ be a complex projective manifold of general type,
with {\em canonical divisor} $K_X$. The {\em canonical ring} 
$R(X,K_X)$
is expected to be finitely generated, and if it is, 
we would obtain a natural birational map 
$$
\Phi\colon X\dashrightarrow Y:=\Proj R(X,K_X).
$$ 
The birational model $Y$ is called the {\em canonical model} 
of $X$. It depends only on the birational class of $X$ and 
it has a canonical polarization, but it has singularities 
in general. For example, $Y$ may have some Du Val 
singularities in dimension two. The singularities that may 
appear on $Y$ were coined {\em canonical singularities}
by Reid~\cite{Reid79}. 

To get to the formal definition of canonical singularities,
let us take a closer look at what $\Phi$ does for $K_X$.
By Hironaka's resolution of singularities, there is 
a {\em Hironaka hut}
\[ \xymatrix{
& X' \ar[dl]_f \ar[dr]^g &  \\
X \ar@{.>}[rr]^\Phi  & & Y
}\]
that is $X'$ is a projective manifold, 
$f,g$ are birational morphisms and $\Phi=g\circ f^{-1}$.
By definition, $K_X$ is the divisor $(\omega)$ of zeros 
and poles of a non-zero top rational differential form 
$\omega\in \wedge^{\dim(X)}\Omega^1_X\otimes_{\mathbb C} 
{\mathbb C}(X)$.
Denote $K_{X'}=(f^*\omega)$ and $K_Y=(g_*f^*\omega)$. 
The latter is a well defined Weil divisor, since $Y$ is 
normal. Since $X$ has no singularities, the divisor 
$A_f=K_{X'}-f^*(K_X)$ is effective and supported 
by the exceptional locus of $f$. Equivalently, the natural
map $f_*\colon R(X',K_{X'})\to R(X,K_X)$ is an isomorphism.
In particular, $g\colon X'\to Y$ is the canonical model 
of $X'$. Since $g$ is a morphism and $K_{X'}$ is a big
divisor, it follows that there exists $r\ge 1$
such that $\vert rK_Y\vert$ defines a projectively
normal embedding, 
and $A_g=\frac{1}{r}(rK_{X'}-g^*(rK_Y))$ is effective
and supported by the exceptional locus of $g$. In
particular, $g_*\colon R(X',K_{X'})\to R(Y,K_Y)$ is 
also an isomorphism:
\[ \xymatrix{
& R(X',K_{X'}) \ar[dl]_\simeq \ar[dr]^\simeq &  \\
R(X,K_X)  &  & R(Y,K_Y)
}\]
Reid~\cite{Reid79} called a normal germ $P\in Y$ 
a canonical singularity if $A_g$ is well defined and effective, 
for a resolution of singularities $g\colon X'\to Y$. 
The coefficients of the $\Q$-divisor $A_g$ are called
{\em discrepancies}.
To understand discrepancies in terms of the 
manifolds that we started with, we go back to our global
setting and note that
$$
K_{X'}=g^*(K_Y)+A_g
$$
is a Zariski decomposition of $K_{X'}$, with positive part
$g^*(K_Y)$ and fixed part $A_g$. Since $\vert rK_Y\vert$ 
defines a linear system free of base points, $rA_g$ 
is the fixed part of the linear system 
$\vert rK_{X'}\vert$. Finally, it turns out that 
$f^*(K_X)=g^*(K_Y)+(A_g-A_f)$ is a Zariski decomposition.


\subsection{Log canonical models of open manifolds}

Let $U$ be a complex quasi-projective manifold of general 
type, in the sense of Iitaka~\cite{Iitaka84}. By 
Hironaka's resolution of singularities, there exists 
an open embedding $U\subset X$ such that $X$ is a proper
manifold, and the complement $X\setminus U=\sum_i E_i$ 
is a divisor with simple normal crossings. The general 
type assumption means that the {\em log canonical divisor}
$K_X+\sum_i E_i$ is big. The {\em log canonical ring} 
$
R(X,K_X+\sum_i E_i)
$ 
is independent of the choice of compactification, and in
fact depends only on the (proper) birational class of $U$.
It is expected to be finitely generated, and if it 
is, we would obtain a natural birational map 
$$
\Phi\colon X\dashrightarrow Y:=\Proj R(X,K+\sum_i E_i).
$$
As before, we can find a Hironaka hut with the extra property
that $\Exc(f)\cup f^{-1}_*(\sum_i E_i)$ is a 
simple normal crossings divisor $\sum_{i'}E_{i'}$.
Denote $B_Y=g_*(\sum_{i'}E_{i'})$. We imitate
the arguments in the compact case, and obtain isomorphisms
\[ \xymatrix{
& R(X',K_{X'}+\sum_{i'}E_{i'}) \ar[dl]_\simeq \ar[dr]^\simeq &  \\
R(X,K_X+\sum_i E_i)  &  & R(Y,K_Y+B_Y)
}\]
Again, $\vert r(K_Y+B_Y)\vert$ defines a normally generated 
embedding for some $r\ge 1$, and we have Zariski 
decompositions $K_{X'}+\sum_{i'}E_{i'}=g^*(K_Y+B_Y)+A_g$ and 
$f^*(K_X+\sum_i E_i)=g^*(K_Y+B_Y)+(A_g-A_f)$. One can 
see that $\Phi^{-1}$ contracts no divisors of $Y$,
and $\Phi_*(\sum_i E_i)=B_Y$. The pair $(Y,B_Y)$ is log 
canonically polarized, and it's singularities are log 
canonical, as we will see shortly. The pair $(Y,B_Y)$ is 
called the {\em log canonical model} of $U$. 

For example, let $U\to M_g$ be a resolution of 
singularities of the moduli space of smooth curves
of genus $g\ge 2$. Then the log canonical model of $U$ is
$(\overline{M}_g,\delta)$, where $\overline{M}_g$ is the moduli
space of stable curves of genus $g$, and $\delta=\overline{M}_g
\setminus M_g$~\cite{CH88}.


\subsection{Log canonical models of log manifolds, log 
discrepancies}

{\em Log manifolds} are the bridge between open and compact 
manifolds. 
They are pairs $(X,\sum_i b_i E_i)$, where $X$ is nonsingular, 
the $E_i$'s are nonsingular divisors intersecting transversely, 
and $b_i\in [0,1]\cap \Q$ for all $i$. We call $\sum_i b_i E_i$
the boundary of the log manifold, and denote it by $B$.
Suppose moreover that $(X,B)$ is of general type, that
is the log canonical divisor $K_X+B$ is big. The log canonical 
ring $R(X,K_X+B)$ is expected to be finitely generated, and 
if it is, we would obtain a birational map
$$
\Phi\colon X \dashrightarrow Y:=\Proj R(X,K_X+B).
$$
Again, we construct a Hironaka hut with the extra property
that $\Exc(f)\cup f^{-1}_*(\sum_i E_i)$ is a 
simple normal crossings divisor. Let $\cup_j F_j$ be the
exceptional locus of $f$ and denote 
$B_Y=g_*(f^{-1}_*B+\sum_j F_j)$. We imitate
the previous argument, and obtain isomorphisms
\[ \xymatrix{
& R(X',K_{X'}+f^{-1}_*B+\sum_j F_j) \ar[dl]_\simeq \ar[dr]^\simeq &  \\
R(X,K_X+B)  &  & R(Y,K_Y+B_Y)
}\]
There exists $r\in \Z_{\ge 1}$ such that $\vert r(K_Y+B_Y)\vert$ 
defines a normally generated embedding, and we have Zariski 
decompositions
\begin{align*}
K_{X'}+f^{-1}_*B+\sum_j F_j & = g^*(K_Y+B_Y)+A_g \\
f^*(K_X+B) & = g^*(K_Y+B_Y)+(A_g-A_f).
\end{align*}
One can see that $\Phi^{-1}$ contracts no divisors of $Y$,
and $\Phi_*(B)=B_Y$. The birational model
$
\Phi\colon (X,B) \dashrightarrow (Y,B_Y)
$
is called the {\em log canonical model} of $(X,B)$. It is 
polarized by the {\em log canonical divisor} $K_Y+B_Y$
(a $\Q$-divisor), 
and its singularities are called {\em log canonical 
singularities}.
The coefficients of the $\Q$-divisor $A_g$ are called 
{\em log discrepancies}. 


\section{Log varieties, minimal log discrepancies}

Log varieties with log canonical 
singularities are locally open subsets of log canonical models. 
For technical purposes, it is better to work in a slightly
more general context, such as non-rational boundaries 
(to take convex hulls and limits of divisors), or even 
non-log canonical singularities (to construct a flat log 
structure at a prescribed point of a polarized manifold, 
cf \S 4-A). 

\begin{defn} 
A {\em log variety} $(X,B)$ is a complex normal 
variety $X$ endowed with an effective $\R$-Weil divisor 
$B=\sum_i b_i E_i$ such that $K_X+B$ is $\R$-Cartier.
\end{defn}

Recall that the canonical divisor $K_X=(\omega)$ is the 
Weil divisor of zeros and poles of a non-zero top rational 
differential form $\omega$ (it depends on the choice of $\omega$, 
but only up to linear equivalence). The $E_i$'s are prime 
divisors and the $b_i$'s are non-negative real numbers. 
The $\R$-Cartier assumption means that locally on $X$, 
$K_X+B$ equals a finite sum $\sum_i r_i (\varphi_i)$, where 
$r_i\in \R$ and $\varphi_i\in {\mathbb C}(X)^\times$.

Let now $\mu\colon X'\to X$ be birational morphism, and 
$E\subset X'$ a prime divisor. We use the same
form to define the canonical class of $X'$, that is
$K_{X'}=(f^*\omega)$. The {\em log discrepancy} of 
$(X,B)$ at $E$ is defined as
$$
a(E;X,B)=\mult_E(K_{X'}+E-\mu^*(K_X+B))\in \R.
$$
The log discrepancy depends only on the valuation that $E$
induces on ${\mathbb C}(X)$. We call such valuations
{\em geometric}, and denote $c_X(E)=\mu(E)$. For example,
if $E$ is a prime divisor in $X$, then $a(E;X,B)=1-\mult_E(B)$.

\begin{defn} The {\em minimal log discrepancy} of a log
variety $(X,B)$ at a Grothendieck point $\eta\in X$ is 
defined as
$$
a(\eta;X,B)=\inf\{a(E;X,B);c_X(E)=\bar{\eta}\}\in \{-\infty\}
\cup \R_{\ge 0}.
$$
The {\em global minimal log discrepancy} is
$
a(X,B)=\inf_{\eta\in X}a(\eta;X,B).
$
\end{defn}

\begin{exmp} $a(0;{\mathbb C}^d)=d$.
\end{exmp}

The reader may check that $a(\eta;X,B)<0$ implies
$a(\eta;X,B)=-\infty$. Otherwise, $a(\eta;X,B)$ is a 
non-negative real number, and the infimum is a minimum.
If $\eta$ has codimension one, then 
$a(\eta;X,B)=1-\mult_\eta(B)$. Otherwise, construct
a log resolution $\mu\colon X'\to X$ such that 
$\mu^{-1}(\bar{\eta})$ is a divisor, and 
$\mu^{-1}(\bar{\eta}), \mu^{-1}_*B$ and 
$\Exc(\mu)=\cup_j F_j$ are all suported by a simple 
normal crossings divisor. Then 
$a(\eta;X,B)=\min_{\mu(F_j)=\eta}a(F_j;X,B)$.
The log pullback formula 
$$
\mu^*(K_X+B)=K_{X'}+\mu^{-1}_*B+\sum_j (1-a(F_j;X,B))F_j
$$
shows that $a(\eta;X,B)\in \frac{1}{r}\Z$ if $r(K_X+B)$ 
is a Cartier divisor near $\eta$.

\begin{defn} A log variety $(X,B)$ has 
{\em log canonical singularities at} $\eta$ if 
$a(\eta;X,B)\ge 0$. We say that $(X,B)$ has 
{\em log canonical singularities} if $a(X,B)\ge 0$.
\end{defn}

\begin{exmp} Consider the log variety 
$(X,\sum_i b_i E_i)$, where $X$ is a manifold,
$\sum_i E_i$ a simple normal crossings divisor 
and $b_i\in [0,1]$ for all $i$. 
Then $a(X,\sum_i b_i E_i)=\min_i (1-b_i)$.
\end{exmp}

\begin{exmp} Let $X$ be a toric variety and 
$X\setminus T=\cup_i E_i$ the complement of the torus.
Then $(X,\sum_i E_i)$ is a log variety with $a(X,\sum_i E_i)=0$
and $K_X+\sum_i E_i=0$.
\end{exmp}

If $a(X,B)\ge 0$, the infimum in its definition is 
also a minimum. With the notation above,
$
a(X,B)=\min(\min_i a(E_i;X,B), \min_j a(F_j;X,B)). 
$ 
Moreover, the formula
$$
K_{X'}+\mu^{-1}_*B+\sum_j F_j=\mu^*(K_X+B)+
\sum_j a(F_j;X,B)F_j
$$
becomes a Zariski decomposition of the relative log manifold 
of general type $(X',(\mu^{-1})_*B+\sum_j F_j)\to X$.
Its log canonical model in the sense of \S 1 is $(X/X,B)$,
where $X/X$ is the identity morphism.

\begin{defn} Let $(X/S,B)$ be a log variety $(X,B)$ with 
$a(X,B)\ge 0$, endowed with a projective morphism 
$\pi\colon X\to S$. We say that $(X/S,B)$ is a 
\begin{itemize}
\item {\em log Fano} if $-(K_X+B)$ is $\pi$-ample.
\item {\em log Calabi-Yau} if $K_X+B$ is $\pi$-numerically trivial.
\item {\em log canonical model} if $K_X+B$ is $\pi$-ample.
\end{itemize}
\end{defn}

It is useful to observe that a log variety $(X,B)$ has log 
canonical singularities at $\eta$ if and only if $(X/X,B)$ 
is a log canonical model in a neighborhood of $\eta$ 
(cf. Example~\ref{lce}). Note however that the three 
geometric types above coincide if $\pi$ is the identity map.

\begin{defn} A germ of log variety $P\in (X,B)$ is called
{\em flat} if $a(P;X,B)=0$.
\end{defn}

The typical example of flat germ is 
$0\in ({\mathbb C}^d,\sum_{i=1}^d H_i)$, where $H_i$
are the coordinate hyperplanes. Our terminology is inspired
by an analogy between germs and projective manifolds, 
where if the minimal log discrepancy corresponds to the 
Kodaira dimension, flat germs correspond to manifolds with 
Kodaira dimension zero (an elliptic curve, for example).
Also, note that $0\in ({\mathbb C},1\cdot 0)$ is the 
only flat log germ in dimension one.


\subsection{Examples of minimal log discrepancies}\label{exam}

Minimal log discrepancies can be easily computed for log
varieties $(X,B)$ such that $X$ is a toric variety and $B$ 
is supported by the complement of the torus (see~\cite{Oda88}
for standard terminology on toric varieties). We only consider 
here $\Q$-factorial, log canonical toric germs of log varieties
$$
P\in (X,B)=(T_N\emb(\sigma),\sum_{i=1}^d b_i H_i).
$$
They are in one-to-one correspondence with the following
data:
\begin{itemize}
\item $\sigma=\{x\in \R^d; x_1,\ldots,x_d\ge 0\}$.
\item $N\subset \R^d$ is a lattice, containing  
$(1,0,\ldots,0),\ldots,(0,\ldots,0,1)$ as primitive vectors.
\item $(b_1,\ldots,b_d)\in [0,1]^d$.
\end{itemize}
The following basic facts provide lots of examples of minimal
log discrepancies:
\begin{itemize}
\item[(a)] $a(\eta_{H_i};X,B)=1-b_i$.
\item[(b)] Let $x\in N^{prim}\cap \sigma$ be a primitive vector. 
Then $x$ defines a barycentric subdivision $\Delta_x$ 
of $\sigma$, and the exceptional locus of the induced 
birational map $T_N\emb(\Delta_x)\to T_N\emb(\sigma)$
is a prime divisor $E_x$. Then 
$a(E_x;X,B)=\sum_{i=1}^d(1-b_i)x_i$.
\item[(c)] Log resolutions exists in the toric category.
Therefore minimal log discrepancies can be computed using 
only valuations $E_x$ as in (b).
\item[(d)] The point $P$ is the unique fixed point of the 
torus action. Its minimal log discrepancy is computed as
follows
$$
a(P;X,B)=\min\{\sum_{i=1}^d(1-b_i)x_i;x\in N\cap \Int(\sigma)\}.
$$
\item[(e)] Let $P\in C\subset X$ be the toric cycle corresponding 
to a face $\tau\prec \sigma$. The minimal log discrepancy at
its generic point is 
$$
a(\eta_C;X,B)=\min\{\sum_{i=1}^d (1-b_i)x_i;x\in N\cap \relint(\tau)\}.
$$
\item[(f)] The global minimal log discrepancy is computed as follows
$$
a(X,B)=\min\{\sum_{i=1}^d (1-b_i)x_i;x\in N\cap \sigma\setminus 0\}.
$$
\item[(g)] In all minimums above, it suffices to consider only the 
finitely many lattice points $x\in N\cap [0,1]^d$.
\end{itemize}

\begin{exmp} Suppose $N=\Z^d$, that is $X={\mathbb C}^d$ and the 
$H_i$'s are the coordinate hyperplanes. For the cycle 
$C:(x_1=\cdots=x_s=0)$, we have $a(\eta_C;X,B)=s-b_1-\cdots-b_s$.
\end{exmp}

\begin{exmp} Suppose $B=0$. Since $\sigma$ is fixed, only the 
lattice $N$ varies.

\begin{itemize}
\item[(i)] Take $N=\Z^2+\Z(\frac{1}{q},\frac{q-1}{q})$, for some
integer $q\ge 2$. The surface germ $P\in X$ is a $A_{q-1}$-singularity.
We compute
$
N\cap (0,1]^2=\{(\frac{k}{q},\frac{q-k}{q});
1\le k\le q-1\}\cup\{(1,1)\},
$
so $a(P;X)=1$.

\item[(ii)] Take $N=\Z^2+\Z(\frac{1}{k},\frac{1}{k})$, for some
positive integer $k$. As above, we compute $a(P;X)=\frac{2}{k}$.

\item[(iii)] Take $N=\Z^3+\Z\frac{1}{q}(1,p,q-p)$, where $p,q$ 
are integers with $1\le p\le q-1,\gcd(p,q)=1$. Then $P\in X$ 
is a terminal $3$-fold singularity, with $a(P;X)=1+\frac{1}{q}$.

\item[(iv)] Take $N=\Z^3+\Z\frac{1}{2q}(1,q,1+q)$, with $q\ge 1$. 
Then $P\in X$ is a $3$-fold singularity with 
$a(P;X)=1+\frac{1}{q}$. 
This germ has the minimal log discrepancy of a terminal singularity, 
but it's not terminal, since it is not an isolated singularity. 
The singular locus of $X$ is $C_2:(x_1=x_3=0)$, and 
$a(\eta_{C_2};X)=\frac{2}{q}$.
\end{itemize}
\end{exmp}

\begin{exmp} For each point $x\in \Q^d\cap (0,1]^d$,
we construct a $d$-dimensional germ of toric log variety 
$P_x\in (X_x,B_x)$, with minimal log discrepancy
$$
a(P_x;X_x,B_x)=\min_{n\ge 0}\sum_{i=1}^d(1+nx_i-\lceil nx_i\rceil).
$$
Choose a positive integer $q$ such that $qx\in \Z^d$, and 
define integers
$$
n_j=\frac{\gcd(q,qx_1,\ldots,\widehat{qx_j},\ldots,qx_d)}{
\gcd(q,qx_1,\ldots,qx_d)} \ (1\le j\le d).
$$
These integers are independent of the choice of $q$.
Let $\sigma\subset \R^d$ be the standard positive cone, 
and $P_x\in X_x:=T_{\Z^d+\Z\cdot x}\emb(\sigma)$ 
the unique point fixed by the torus.
The primitive vectors of the lattice $\Z^d+\Z\cdot x$ along 
the rays of $\sigma$ are 
$$
\frac{1}{n_1}(1,0,\ldots,0),\ldots,\frac{1}{n_d}(0,\ldots,0,1).
$$ 
They define invariant
prime divisors $H_1,\ldots,H_d$ in $X$. Set
$B_x=\sum_{i=1}^d (1-\frac{1}{n_i})H_i$. 

It is natural to ``compactify'' $\Q^d\cap (0,1]^d$ to 
$\Q^d\cap [0,1]^d\setminus 0$, in order to study the 
limiting behaviour of $a(P_x;X_x,B_x)$~\cite{Bor97,A06}. 
\end{exmp}

The minimal log discrepancy of a toric germ 
is the ``local'' version of Minkowski's first minimum 
of a convex body about the origin (see~\cite{Cassels57}
for the definition and basic properties of the latter 
invariant). In our setting, if we denote by $\Delta$
the convex set 
$\{x\in \R^d; x_1,\ldots,x_d\ge 0,\sum_{i=1}^d (1-b_i)x_i\le 1\}$, 
then 
$$
a(P;X,B)=\sup\{t\ge 0;N\cap \Int(t\Delta)=\emptyset\}.
$$
See~\cite{Bor00} for more on 
minimal log discrepancies versus lattice-point-free 
convex bodies.


\section{Problems on minimal log discrepancies}


Minimal log discrepancies originate in the problem of 
termination of log flips: starting with a given log variety, 
can we perform log flips infinitely many times? 
Log flips are surgery operations which preserve codimension one cycles, 
and improve the singularities of higher codimensional cycles. 
As a measure of this improvement, log discrepancies may only increase
after a log flip, and some of them increase strictly~\cite{Sh85}. 
This is the heuristic behind the termination 
of a sequence of log flips, and it lead Shokurov~\cite{Sh85,Sh88,Sh96} 
to question the existence of an infinite increasing sequence of 
minimal log discrepancies. 

First, we fix a log variety $(X,B)$, and investigate the set of minimal 
log discrepancies of all prime cycles of $X$~\cite{Ab99}. 
The formula $a(\eta_C;X,B)=a(P;X,B)-\dim(C)$, for 
a general closed point $P$ of a given prime cycle $C\subset X$, 
shows that closed points contain the essential information. 
Consider now the minimal log discrepancy 
$a(P;X,B)$ as a function on the set of the closed points $P\in X$. This 
function has a finite image, and in particular the set of minimal log 
discrepancies of all prime cycles of $X$ is {\em finite}. Moreover, the
level sets $\{P\in X;a(P;X,B)\le t\} \ (t\ge 0)$ are constructible. 
Simple examples, such as a Du Val singularity $P\in X$, with $a(x;X)=2$ 
for $x\ne P$, and $a(P;X)=1$, suggest that these level sets are in fact 
closed.

\begin{conj}[\cite{Aa99}]\label{lsc} 
The minimal log discrepancy $a(P;X,B)$ is 
lower semi-continuous as a function on the closed points $P$ of $X$.
\end{conj}

This behaviour is confirmed in several cases: 
a) $\dim(X)\le 3$~\cite{Aa99,Ab99}; b) $(X,B)$ is a toric log 
variety~\cite{Ab99}; c) $X$ is a local complete intersection~\cite{EM,EMY}.
Also, it is equivalent to the inequality $a(P;X,B)\le a(\eta_C;X,B)+1$,
for every closed point on a curve in $X$~\cite{Ab99}.

Now consider the general case, when log flips change the log variety 
$(X,B)$ in codimension at least two. The coefficients of the boundary are 
preserved, so we may assume that they belong to a given finite set. 
More generally, let $\cB\subset [0,1]$ be a set satisfying the 
{\em descending chain condition} ($\cB=\{1-\frac{1}{n};n\ge 1\}\cup\{1\}$ 
is the typical example), and define
$$
\Mld(d,\cB)=\{a(P;X,B); \dim(X)=d, 
\textrm{ coefficients of }B \textrm{ belong to } \cB\}.
$$
The set $\Mld(1,\cB)=\{1-b;b\in \cB\}$ clearly satisfies the 
ascending chain condition.

\begin{conj}[Shokurov~\cite{Sh88,Sh96}]\label{MC}
The following properties hold:
\begin{itemize}
\item[(1)] $\Mld(d,\cB)$ satisfies the ascending chain condition.
\item[(2)] $a(P;X,B)\le \dim(X)$. Moreover, if $a(P;X,B)>\dim(X)-1$, 
then $P\in X$ is a nonsingular point and $a(P;X,B)=\dim(X)-\mult_P(B)$.
\item[(3)] Assume that $\cB\cap [0,1-\frac{1}{k}]$ is a finite set for 
every $k\ge 1$. Then the accumulation points of $\Mld(d,\cB)$ are 
included in $\Mld(d-1,\cB')$, for a suitable set $\cB'$.
\end{itemize}
\end{conj}
This conjecture was confirmed for surfaces~\cite{Sh91,Ale93}, and 
toric log varieties~\cite{Bor97,A06}. By the classification of 
terminal $3$-fold singularities, 
$\Mld(3,\{0\})\cap (1,+\infty)=\{1+\frac{1}{q};q\ge 1\}\cup\{3\}$
~\cite{Kaw93, Mrk96}.
Also, (2) holds if $X$ is a local complete intersection~\cite{EM,EMY}.
Shokurov~\cite{Sh04} reduced the global problem of the termination 
of flips to the two local problems~\ref{lsc} and ~\ref{MC}.(1).
See~\cite{BS06} for more on this.

An interesting problem posed by Shokurov is to relate
the minimal log discrepancy and the Cartier index of 
a singularity. Suppose $P\in X$ is the germ of a $d$-fold 
with $a(P;X)\ge 0$. If $nK_X\sim 0$ and Conjecture~\ref{MC}.(2) 
holds, then the minimal log discrepancy can take at most 
finitely many values:
$
a(P;X)\in \{0,\frac{1}{n},\ldots,\frac{nd}{n}\}.
$
Conversely, is there an integer $n$, depending only 
on $d$ and $a(P;X)$, such that $nK_X\sim 0$? The answer is
positive if $d=2$ (Shokurov, unpublished). An important
special case of this problem is $a(P;X)=0$. If $d=2$, then 
$n\in \{1,2,3,4,6\}$~\cite{Sh92}. If $d=3$, then $\varphi(n)\le 
20$ and $n\ne 60$, where $\varphi$ is the Euler 
number~\cite{Ishii98}. See also~\cite{Fuj01} for a higher 
dimensional reduction to a global problem on log
Calabi-Yau varieties in one dimension less. The boldest
conjecture here is the following.

\begin{conj}[Shokurov~\cite{Sh00}] Let $P\in (X,B)$ be a 
log germ such that $a(P;X,B)=0$ and $B$ has coefficients 
in a set $\cB\subset [0,1]\cap \Q$ satisfying the descending
chain condition. Then $n(K+B)\sim 0$ for some positive 
integer $n$, depending only on $\dim(X)$ and $\cB$. 
\end{conj}

A useful formula in inductive arguments in the log 
category is a comparison of minimal log discrepancies under 
adjunction, called precise inversion of adjunction.

\begin{conj}[Shokurov~\cite{Sh92}, Koll\'ar~\cite{Kol92}] Let
$P\in S\subset (X,B)$ be the germ of a log variety and a normal
prime divisor $S$ with $\mult_S(B)=1$. By adjunction, we have
$(K_X+B)\vert_S=K_S+B_S$. Then $a(P;X,B)=a(P;S,B_S)$.
\end{conj}

It follows from the Log Minimal Model Program if 
$a(P;X,B)\le 1$~\cite{Kol92}, and it holds if $X$ is 
a local complete intersection~\cite{EM,EMY}.

Minimal log discrepancies appear naturally in global
contexts, such as Fujita's Conjecture on adjoint linear
systems or the boundedness of log Fano varieties.

\begin{conj}[A. and L. Borisov~\cite{BB}-Alexeev~\cite{Ale94}]\label{ba}
Let $\epsilon\in (0,1]$ and $d\in \Z_{\ge 1}$. Then log 
varieties $X$ with $-K_X$ ample, $a(X)\ge \epsilon$ and
$\dim(X)=d$, form a bounded family.
\end{conj}

This is known in several cases: a) $X$ is toric~\cite{BB};
b) $X$ nonsingular~\cite{KMM92}; c) $d=2$~\cite{Ale94}; d) $d=3$,
$\epsilon=1$~\cite{Kaw92,KMMT}; e) $d=3$, and the index of $K_X$ is 
fixed~\cite{Bor01}.


\subsection{Toric case}

In the assumptions and notations of \S~\ref{exam}, we 
illustrate some of the local problems on minimal log 
discrepancies.
For {\em lower semi-continuity}, it is enough to see that
$a(P;X,B)\le a(\eta_C;X,B)+1$ for a torus-invariant curve $P\in C$.
Suppose $C$ corresponds to the face $\tau=\sigma\cap (x_d=0)$.
There exists $(x',0)\in N^{prim}\cap \relint(\sigma)$ such that 
$a(\eta_C;X,B)=\sum_{i=1}^{d-1}(1-b_i)x_i$. Then 
$(x',1)\in N\cap \Int(\sigma)$ and there exists 
$x\in N^{prim}$ and a positive integer $l$ with
$lx=(x',1)$. We have
$$
a(E_x;X,B)\le la(E_x;X,B)=\sum_{i=1}^{d-1}(1-b_i)x'_i+1-b_d 
\le a(\eta_C;X,B)+1.
$$
Therefore $a(P;X,B)\le a(\eta_C;X,B)+1$.

For {\em precise inversion of adjunction}, suppose 
$B=\sum_{i=1}^{d-1}b_i H_i+H_d$. 
Then $S=H_d$ is the toric variety $T_{N_d}\emb(\sigma_d)$, where
$\sigma_d=\{x\in \R^{d-1};x_1,\ldots,x_{d-1}\ge 0\}$ and
$N_d=\{x\in \R^{d-1}; ^\exists t\in \R, (x,t)\in N\}$. To bring
this to the normal form in \S~\ref{exam}, note that there are
positive integers $n_1,\ldots,n_{d-1}$ such that 
$\frac{1}{n_i}(0,\ldots,\stackrel{i}{1},\ldots,0)$ are primitive
vectors of $N_d^{prim}$. Then $S=T_{N'}\emb(\sigma')$, where 
$N'=\{x'\in \R^{d-1};(n_1x'_1,\ldots,n_{d-1}x'_{d-1})\in N_d\}$
and $\sigma'$ is the usual positive cone. Let $H'_1,\ldots,
H'_{d-1}$ the torus invariant prime divisors of $S$. The key observation
is that the log canonical divisor $K_X+B=\sum_{i=1}^{d-1}-(1-b_i)H_i$ is 
independent of $H_d$. The boundary induced on $S$
by adjunction is $B_S=\sum_{i=1}^{d-1}(1-\frac{1-b_i}{n_i})H'_i$, and 
the equality $a(P;X,B)=a(P;S,B_S)$ is clear.

Finally, for the {\em ascending chain condition}, assume by contradiction 
that we have a strictly increasing sequence
$
a^1<a^2<a^3<\cdots,
$
where $a^n=a(P^n;T_{N^n}\emb(\sigma))$ for $n\ge 1$. For simplicity, 
assume that the boundary is zero, so only the lattice
changes. We may find $x^n\in (0,1]^d\cap N^n$ such that 
$a^n=\sum_{i=1}^d x_i^n$. In particular, $a^n\le d$ for all $n$. 
Consider now the strictly increasing sequence of open sets
$
U^n=\{x\in (0,+\infty)^d; \sum_{i=1}^d x_i<a^n\}.
$
By~\cite{Lawrence91}, $G^n=\{x\in \R^d; U^n\cap (\Z^d+\Z x)=\emptyset\}$
is the union of finitely many closed subgroups containing $\Z^d$
(the Flatness Theorem of Khinchin~\cite{KL} gives an alternative proof). 
We have $G^n\supsetneq G^{n+1}$ since 
$a^n<a^{n+1}$ and $x^n\in G^n\setminus G^{n+1}$, so we obtain 
a strictly decreasing sequence of finite unions of closed subgroups 
containing $\Z^d$. This is impossible, since the latter set 
satisfies the descending chain condition.


\subsection{Methods}


The toric case (see~\cite{MMM,Bor99}) suggests that behind 
the ascending chain condition of minimal log discrepancies 
lies a deeper fact, the boundedness of singularities 
with minimal log discrepancy bounded away from zero. 
Some log canonical singularities are classified 
in low dimension, but in general we could only expect general structure 
theorems and boundedness results in terms of minimal log discrepancies.
For example, Du Val singularities are classified as follows:
$A_n,D_n,E_6,E_7,E_8$. From the above point of view, Du Val
singularities are nothing but surface singularities having minimal log
discrepancy at least $1$, and they come in two types: a $1$-dimensional
series with two components ($A$ and $D$), and a $0$-dimensional 
series ($E$). 

The known method for bounding germs $P\in X$ is to study the 
singularities at $P$ of the linear systems $\vert -mK_X\vert\ (m\ge 1)$,
and reduce this local problem to the global problem of bounding 
log Fano or log Calabi-Yau varieties in one dimension 
less~\cite{Sh00,Pro01}. Given that \S 1 suggests that minimal log 
discrepancies are invariants of objects of general type, 
it also natural to investigate the singularities at $P$ of 
the linear systems $\vert mK_X\vert \ (m\ge 1)$ and 
$\vert {\mathfrak m}_{P,X}^m\vert \ (m\ge 1)$, and 
relate germs with log canonical models in one dimension less.

Also, it is possible that minimal log discrepancies can be understood 
from several points of view: analytic, birational, motivic or p-adic. The 
motivic interpretation of minimal log discrepancies is known in the 
case when the canonical divisor is $\Q$-Cartier~\cite{Mircea02,Yasuda03}. 
As for the
analytic side, the description of log discrepancies as the coefficients
of a Zariski decomposition suggests an interpretation of minimal log
discrepancies in terms of Lelong numbers. For example, the upper bound 
of Conjecture~\ref{MC}.(2) is essentially equivalent to the following
problem: let $(X,B)$ be a log manifold of general type
having a Zariski decomposition $K_X+B=P+F$, such that $F$ is 
supported by the components of $B$ with coefficient one. 
Then some coefficient of $F$ is at most $\dim(X)$.


\section{The log canonical threshold}


\begin{defn} Let $(X,B)$ be a log variety, $\eta\in X$
a Grothendieck point such that $a(\eta;X,B)\ge 0$, and $D$ 
an effective $\R$-Cartier $\R$-divisor passing through $\eta$.
The {\em log canonical threshold} of $D$ with respect to
$\eta\in (X,B)$ is defined as
$$
\lct(\eta\in (X,B);D)=\sup\{t\ge 0;a(\eta;X,B+tD)\ge 0\}.
$$
The {\em global log canonical threshold} is  
$
\lct(X,B;D)=\min_{\eta\in X}\lct(\eta\in (X,B);D).
$
\end{defn}

Note that $a(\eta;X,B)\ge 0$ if and only if $(X,B)$ has log
canonical singularities at $\eta$. Since the support of 
$B+tD$ is a fixed divisor $S$, we may compute the log canonical 
threshold on a log resolution of $X$ and $S$. In particular,
$\lct(\eta\in (X,B);D)\in \Q$ if $B$ and $D$ are rational.
The inequality $\lct(P\in (X,B);D)\le \lct(\eta_C\in (X,B);D)$ 
holds for a closed point of a prime cycle $P\in C\subset X$, 
with equality if $P$ is general in $C$.

\begin{exmp} Let $\dim(X)=1$ and $P\in X$. 
Then $\lct(P\in (X,B);D)=\frac{a(P;X,B)}{\mult_P(D)}$.
\end{exmp}

\begin{exmp}[\cite{Demailly01,Hwang05}] 
Let $0\in H\colon (f=0)\subset {\mathbb C}^d$ be a 
hypersurface. Then
$$
\lct(0\in {\mathbb C}^d;H)=
\sup\{t\ge 0; \vert f\vert^{-t}\mbox{ is }L^2 \mbox{ near }0\}.
$$
The reciprocal number $\mu(0;H)=1/\lct(0\in {\mathbb C}^d;H)$,
called Arnold multiplicity, satisfies the 
inequalities $\mu(0;H)\le \mult(0;H)\le d\cdot \mu(0;H)$.
\end{exmp}

\begin{exmp} Let $H_1,\ldots,H_d$ be the coordinate 
hyperplanes in ${\mathbb C}^d$, $b_1,\ldots,b_d\in [0,1]$
and $n_1,\ldots,n_d\in \Z_{\ge 1}$. Then 
\begin{itemize}
\item[(i)] $\lct(0\in ({\mathbb C}^d,\sum_{i=1}^d b_i H_i);
(x_1^{n_1}\cdots x_d^{n_d}))
=\min_{i=1}^d\frac{1-b_i}{n_i}$.
\item[(ii)] $\lct(0\in ({\mathbb C}^d,\sum_{i=1}^d b_i H_i);
(x_1^{n_1}+\cdots +x_d^{n_d}))
=\min(1,\sum_{i=1}^d\frac{1-b_i}{n_i})$.
\end{itemize}
\end{exmp}

\begin{exmp} Log canonical thresholds of non-degenerate
hypersurfaces in toric varieties have a combinatorial
description. Let 
$P\in (X,B)=(T_N\emb(\sigma),\sum_{i=1}^d b_i H_i)$ 
be a toric germ as in \S 2-A, and 
$$
f=\sum_\alpha \lambda_\alpha \chi^{m_\alpha}\in 
{\mathfrak m}_{P,X}\setminus 0.
$$
Here $m_\alpha$ are finitely many points in 
$M\cap \sigma^\vee$, where $M$ is the lattice dual
to $N$ and $\sigma\subset M_\R$ is the cone dual
to $\sigma$. Note that $(1-b_1,\ldots,1-b_d)\in 
\sigma^\vee$. Let $\square$ be the convex hull of 
the $\cup_\alpha(m_\alpha+\sigma^\vee)$. The 
ray $\R_{\ge 0}(1-b_1,\ldots,1-b_d)$ intersects 
$\square$ for the first time in a point
$\mu(1-b_1,\ldots,1-b_d)$. 
We have the following inequality
$$
\lct(P\in (X,B);(f=0))\le \min(1,\frac{1}{\mu}),
$$
and equality holds if the coefficients $\lambda_\alpha$ 
are sufficiently general.
The generality condition can be made explicit, as 
in~\cite{AGV88}.
\end{exmp}

\begin{exmp}\label{lce} The log canonical threshold 
$\lct(P\in (X,B);D)$ is the largest $t\ge 0$ such that 
$(X/X,B+tD)$ is a relative log canonical model near $P$. 
For example, consider the germ of the cuspidal curve 
$0\in C\colon (x^2-y^3=0)\subset {\mathbb C}^2$,
with log canonical threshold $\frac{5}{6}$.
The cusp can be resolved by a composition of three
blow-ups $\mu\colon X\to {\mathbb C}^2$. Let $E_1,E_2,E_3$ be
the proper transforms on $X$ of the exceptional divisors, 
in the order of their appearance. We obtain a singular model 
$X\to Y\to {\mathbb C}^2$ by contracting $E_1$ and $E_2$,
and let $C'$ and $E'_3$ be the proper transforms of $C$ 
and $E_3$ on $Y$. The relative log variety 
$
(X,E_1+E_2+E_3+t\mu^{-1}_*C)\to {\mathbb C}^2
$
is of general type for $t\in [0,1]\cap \Q$, and 
its log canonical model is 
$({\mathbb C}^2,tC)\to {\mathbb C}^2$ for $t\le \frac{5}{6}$, 
and $(Y,E'_3+tC')\to {\mathbb C}^2$ for $t>\frac{5}{6}$.
On the other hand, $({\mathbb C}^2,tC)\to {\mathbb C}^2$
is a log canonical model on ${\mathbb C}^2\setminus 0$, 
for every $t\in [0,1]\cap \Q$.
\end{exmp}


\subsection{Problems on log canonical thresholds}


One can classify a singularity by constructing a 
{\em flat log structure} on it, in an effective way. Log 
canonical thresholds appear naturally in this construction, 
as the coefficients of the boundary of the flat log structure. 

For example, consider a log germ $P\in (X,B)$ with 
$a(P;X,B)\ge 0$. If $a(P;X,B)=0$, the log germ is flat.
Otherwise, for a general hypersurface 
$D_1\in \vert {\mathfrak m}_{P,X}\vert$, the log canonical 
threshold 
$\gamma_1=\lct(P\in (X,B);D_1)$ is defined in such a way 
so that $a(P;X,B+\gamma D_1)\ge 0$, and there exists a 
minimal cycle $C\ni P$ such that $a(\eta_C;X,B+\gamma D)=0$. 
If $C=\{P\}$, then $P\in (X,B+\gamma D)$ is flat. Otherwise, 
we repeat this process (at most 
$\dim(C)$ times) until we obtain a flat log structure 
$0\in (X,B+\sum_i \gamma_i D_i)$. 
Flat germs with certain restrictions on their boundaries 
are expected to be bounded in a certain sense (cf. \S 3), 
so for effective results we should also control
the $\gamma_i$'s. Shokurov observed 
that these coefficients satisfy the ascending chain condition 
in dimension two, and used this to construct $3$-fold 
log flips.

\begin{conj}[Shokurov~\cite{Sh92}, Koll\'ar~\cite{Kollar97}]\label{lc}
For a positive integer $d$ and a set $0\in \cB\subset [0,1]$ 
satisfying the descending chain condition, define
$$
\Lct(d,\cB)=\{\lct(P\in (X,B);D);\dim(X)=d, 
\textrm{coefficients of }B \textrm{ belong to } \cB,
D\in \vert {\mathfrak m}_{P,X}\vert\}.
$$
The following properties hold:
\begin{itemize}
\item[(1)] The set $\Lct(d,\cB)$ satisfies the ascending chain 
condition. 
\item[(2)] Assume that $\cB\cap [1,1-\frac{1}{k}]$ is a finite
set for every $k\ge 1$. Then the set of accumulation points
of $\Lct(d,\cB)$ is $\Lct(d-1,\cB')\setminus \{1\}$, for a 
suitable set $\cB'\subset [0,1]$.
\end{itemize}
\end{conj}

This limiting behaviour is known in several cases:
a) $d=2$~\cite{Sh92}; b) $d=3$~\cite{Ale94,MP04}; 
c) $P\in X$ is toric, $B$ is invariant and 
$D\in \vert {\mathfrak m}_{P,X}\vert$ is a general 
hypersurface~\cite{Ishii}. The set $\Lct_d:=\Lct(d,\{0\})$ 
has an explicit description in some cases:
a) $\Lct_1=\{\frac{1}{n};n\in \Z_{\ge 1}\}$;
b) $\Lct_2\cap [\frac{1}{2},1]=\{\frac{1}{2}+\frac{1}{n};n\ge 3\}
\cup \{1\}$~\cite{Kuwata99}; c)
$\Lct_3\cap [\frac{41}{42},1]=\{\frac{41}{42},1\}$~\cite{Kollar94}.
Conjecture~\ref{lc} would imply that the number
$1-\epsilon_d=
\max\Lct_d\cap (0,1)$ is well defined, and Koll\'ar~\cite{Kollar94}
suggests that $\epsilon_d$ is the minimal degree
of a log canonical model $(Y,B_Y)$ of dimension $d-1$, with 
$B_Y$ having coefficients in $\{1-\frac{1}{n};n\ge 1\}\cup\{1\}$.
For example, $\epsilon_1=\frac{1}{2},\epsilon_2=\frac{1}{6},
\epsilon_3=\frac{1}{42}$.
For the relationship between Conjectures~\ref{MC},\ref{ba} 
and ~\ref{lc}, see~\cite{BS06}.

Flat log structures have global applications as well.
Given a closed point $P$ on a polarized manifold $(X,H)$, 
does there exists $n\ge 1$
and $D\in \vert nH\vert$ with $a(P;X,\frac{1}{n}D)=0$? 
If such a flat structure exists, then $P$ is not in the base 
locus of any divisor
$L$ such that $L-K_X-H$ ample. This is a powerful technique to
study adjoint line bundles, parallel
to the L$^2$-methods for singular hermitian metrics in 
complex geometry (see~\cite{KMM, Ein97, Kollar97, Aa99} 
and~\cite{Siu95,Demailly01} for the algebraic and analytic 
side of the story, respectively).
As above, log canonical thresholds appear naturally.
For a log variety $(X,B)$, an ample $\Q$-divisor 
$H$ and a Grothendieck point $\eta\in X$, define
$$
\gamma(\eta\in (X,B);H)=
\inf\{n\lct(\eta\in (X,B);D);\ D\in \vert nH\vert,n\ge 1\}
$$
and
$\gamma(X,B;H)=\inf_{\eta\in X}\gamma(\eta \in (X,B);H)$ 
(cf.~\cite{Viehweg95}).

\begin{exmp} Suppose $\dim(X)=1$. For every closed point
$P\in X$ we have 
$$
\gamma(\eta\in (X,B);H)=\frac{a(P;X,B)}{\deg(H)}\le 
\frac{1}{\deg(H)},
$$
and $\gamma(X,B;H)=\frac{a(X,B)}{\deg(H)}$.
\end{exmp}

\begin{lem}[\cite{Sh85,Kollar97}]\label{bi} 
$
\gamma(P\in (X,B);H)\sqrt[d]{H^d}\le d
$
for every closed point $P\in X^d$.
\end{lem}

For Fujita's Conjecture on adjoint linear systems, 
$\gamma(P\in (X,B);H)$ has to be small.
This is achieved by Lemma~\ref{bi}, if $\sqrt[d]{H^d}$ 
is large. 
On the other hand, $\gamma(P\in (X,B);H)$ cannot be 
too small in a bounded context, and this can
be used to bound $\sqrt[d]{H^d}$ from above. 
Some known cases are:
a) $\gamma(\eta_X\in X;-K_X)\ge \frac{1}{d+1}$ if $X$ is
a Fano manifold of dimension $d$~\cite{Hwang05};
b) $\gamma(G;H)\ge 1$ if $\cO_G(H)$ is the Pl\"ucker line 
bundle on the Grassmanian~\cite{Hwang06};
c) $\gamma(A;\Theta)\ge 1$ for a principally polarized
abelian variety~\cite{Kollar95}. 
Pukhlikov~\cite{Pukh04} suggests that a Fano variety $X$ 
with large $\gamma(X;-K_X)$ is birationally rigid.
Hwang~\cite{Hwang05} suggests that for a log canonical model $X$,
$\gamma(\eta_X\in X;K_X)\sqrt[d]{K_X^d}$ is also bounded 
away from zero, only in terms of dimension. The level sets 
and critical points of the functional $\gamma(\eta\in (X,B);H)$, 
and an extension of Faltings' product 
theorem~\cite{Faltings91, Nadel91,Hwang05},
should play a key role in constructing flat log structures.



\begin{thebibliography}{K87}


\bibitem{Ale93}
Alexeev, V.,
{\em Two two-dimensional terminations.}
{Duke Math. J. 69 (1993), no. 3, 527--545.} 

\bibitem{Ale94}
Alexeev, V.,
{\em Boundedness and $K\sp 2$ for log surfaces.}
{Internat. J. Math. 5 (1994), no. 6, 779--810.} 

\bibitem{Aa99}
Ambro, F.,
{\em The Adjunction Conjecture and its applications,}
{preprint math.AG/9903060.}

\bibitem{Ab99}
Ambro, F.,
{\em On minimal log discrepancies,}
{Math. Res. Letters {\bf 6}, (1999) 573--580.}

\bibitem{A06}
Ambro, F.,
{\em The set of toric minimal log discrepancies.}  
{Cent. Eur. J. Math. 4 (2006), no. 3, 358--370.} 

\bibitem{AGV88}
Arnold, V. I.; Gusen-Zade, S. M.; Varchenko, A. N.,
{\em Singularities of differentiable maps. Vol. II. Monodromy 
and asymptotics of integrals. 
Monographs in Mathematics, 83. Birkh\"auser Boston, Inc., Boston, 
MA, 1988.}

\bibitem{BS06}
Birkar, C.; Shokurov, V.,V.,
{\em Mld's vs thresholds and flips.}
{preprint math.AG/0609539.}

\bibitem{BB}
Borisov, A. A.; Borisov, L. A.,
{\em Singular toric Fano three-folds.} 
{(Russian)  Mat. Sb.  183  (1992), no. 2, 134--141;  
translation in  Russian Acad. Sci. Sb. Math.  75  (1993),  
no. 1, 277--283.}

\bibitem{Bor97}
Borisov, A.,
{\em Minimal discrepancies of toric singularities.}  
{Manuscripta Math. 92 (1997), no. 1, 33--45.}

\bibitem{Bor99}
Borisov, A.,
{\em On classification of toric singularities.}
{Algebraic geometry, 9 J. Math. Sci. (New York)
94 (1999), no 1, 1111--1113.}

\bibitem{Bor00}	
Borisov, A.,
{\em Convex lattice polytopes and cones with few 
lattice points inside, from a birational geometry 
viewpoint.} {preprint math.AG/0001109.}

\bibitem{Bor01}
Borisov, A.,
{\em Boundedness of Fano threefolds with log-terminal 
singularities of given index.}  
{J. Math. Sci. Univ. Tokyo 8 (2001), no. 2, 329--342.}

\bibitem{Cassels57}
Cassels, J. W. S.,
{\em An introduction to Diophantine approximation.}
{Cambridge Tracts in Mathematics and Mathematical Physics, 
No. 45. Cambridge University Press, New York, 1957.}

\bibitem{CH88}
Cornalba, M.; Harris, J.,
{\em Divisor classes associated to families of stable 
varieties, with applications to the moduli space of curves.}
{Ann. Sci. \'Ecole Norm. Sup. (4) 21 (1988), no. 3, 455--475.} 

\bibitem{Cutkosky86}
Cutkosky, S, D.,
{\em Zariski decomposition of divisors on algebraic varieties.}
{Duke Math. J. 53 (1986), no. 1, 149--156.} 

\bibitem{Demailly01}
Demailly, J.-P.,
{\em Multiplier ideal sheaves and analytic methods 
in algebraic geometry.} 
{School on Vanishing Theorems and Effective Results 
in Algebraic Geometry (Trieste, 2000),  1--148, ICTP 
Lect. Notes, 6, Abdus Salam Int. Cent. Theoret. Phys., 
Trieste, 2001.}

\bibitem{Ein97}
Ein, L.,
{\em Multiplier ideals, vanishing theorems and applications.} 
{Algebraic geometry---Santa Cruz 1995, 203--219,
Proc. Sympos. Pure Math., 62, Part 1, 
Amer. Math. Soc., Providence, RI, 1997.} 

\bibitem{EMY}
Ein, L.; Musta\c{t}\u{a}, M.; Yasuda, T.,
{\em Jet schemes, log discrepancies and inversion of adjunction.}
{Invent. Math. 153 (2003), no. 3, 519--535.} 

\bibitem{EM}
Ein, L.; Musta\c{t}\v{a}, M.,
{\em Inversion of adjunction for local complete intersection 
varieties.}
{Amer. J. Math. 126 (2004), no. 6, 1355--1365.} 

\bibitem{Faltings91}
Faltings, G.,
{\em Diophantine approximation on abelian varieties.}  
{Ann. of Math. (2)  133  (1991),  no. 3, 549--576.}

\bibitem{Fuj01}
Fujino, O.,
{\em The indices of log canonical singularities.}  
{Amer. J. Math. 123 (2001), no. 2, 229--253.}

\bibitem{Fujita86}
Fujita, T.,
{\em Zariski decomposition and canonical rings of 
elliptic threefolds.}  
{J. Math. Soc. Japan 38 (1986), no. 1, 19--37.}

\bibitem{Hwang05}
Hwang, J.-M.,
{\em Arnold multiplicity of divisors on rational 
homogeneous spaces.}
{Proceedings of Algebraic Geometry Symposium (Kinosaki
2005), A. Ishii (Ed.), pp. $78-81$}.

\bibitem{Hwang06}
Hwang, J.-M.,
{\em Log canonical thresholds of divisors on 
Grassmannians.} {Math. Ann. 334 (2006), no. 2, 413--418.}

\bibitem{Iitaka84}
Iitaka, S.,
{\em Birational geometry of algebraic varieties.}  
{Proceedings of the International Congress of Mathematicians, 
Vol. 1, 2 (Warsaw, 1983), 727--732, PWN, Warsaw, 1984.}

\bibitem{Ishii98}
Ishii, S.,
{\em The quotients of log-canonical singularities by 
finite groups.} 
{Singularities---Sapporo 1998, 135--161,
Adv. Stud. Pure Math., 29, Kinokuniya, Tokyo, 2000.}

\bibitem{Ishii}
Ishii, S.,
{\em Log canonical thresholds for pairs ${\mathbb C}^n$
and non-degenerate hypersurfaces.}{preprint 2002.}

\bibitem{KL}
Kannan, R.; Lov\'asz, L.,
{\em Covering minima and lattice-point-free convex bodies.}
{Ann. of Math. (2) 128 (1988), no. 3, 577--602.} 

\bibitem{Kawamata85}
Kawamata, Y.,
{\em The Zariski decomposition of log-canonical divisors.} 
{Algebraic geometry, Bowdoin, 1985 (Brunswick, Maine, 1985), 
425--433, Proc. Sympos. Pure Math., 46, Part 1, Amer. Math. 
Soc., Providence, RI, 1987.}

\bibitem{KMM} 
Kawamata Y., Matsuda K., Matsuki K., 
{\em Introduction to the minimal model program}, 
{Algebraic Geometry, Sendai, Advanced Studies in Pure Math. 
{\bf 10} (1987), 283--360.}

\bibitem{Kaw92}
Kawamata, Y.,
{\em Boundedness of $\mathbb Q$-Fano threefolds.} 
{Proceedings of the International Conference on Algebra, 
Part 3 (Novosibirsk, 1989), 439--445, Contemp. Math., 131, 
Part 3, Amer. Math. Soc., Providence, RI, 1992.}

\bibitem{Kaw93} 
Kawamata, Y., 
{\em The minimal discrepancy of a $3$-fold terminal 
singularity}, {An appendix to \cite{Sh92}}.

\bibitem{Kol92}
Koll\'ar, J.,
{\em Adjunction and discrepancies.}
{In {\em Flips and abundance for algebraic threefolds}
(Koll\'ar, J Ed.), Ast\'erisque 211, 1992.}

\bibitem{KMM92}
Koll\'ar, J.; Miyaoka, Y.; Mori, S.,
{\em Rational connectedness and boundedness of Fano manifolds.}
{J. Differential Geom. 36 (1992), no. 3, 765--779.} 

\bibitem{Kollar94}
Koll\'ar, J.,
{\em Log surfaces of general type; some conjectures.}  
{Classification of algebraic varieties (L'Aquila, 1992),  
261--275, Contemp. Math., 162, Amer. Math. Soc., 
Providence, RI, 1994.}

\bibitem{Kollar95}
Koll\'ar, J.,
{\em Shafarevich Maps and Automorphic Forms,} 
{Princeton Univ. Press, 1995.}

\bibitem{Kollar97}
Koll\'ar, J\'anos,
{\em Singularities of pairs.} 
{Algebraic geometry---Santa Cruz 1995, 221--287,
Proc. Sympos. Pure Math., 62, Part 1, Amer. Math. Soc., 
Providence, RI, 1997.}

\bibitem{KMMT}
Koll\'ar, J.; Miyaoka, Y.; Mori, S.; Takagi, H.,
{\em Boundedness of canonical $\mathbb Q$-Fano 3-folds.}
{Proc. Japan Acad. Ser. A Math. Sci. 76 (2000), no. 5, 73--77.}

\bibitem{Kuwata99}
Kuwata, T.,
{\em On log canonical thresholds of reducible 
plane curves.}  
{Amer. J. Math. 121 (1999), no. 4, 701--721.}

\bibitem{Lawrence91}
Lawrence, J., 
{\em Finite unions of closed subgroups of the $n$-dimensional 
torus.}
{Applied geometry and discrete mathematics, 433--441, 
DIMACS Ser. Discrete Math. Theoret. Comput. Sci., 4, 
Amer. Math. Soc., Providence, RI, 1991.}

\bibitem{Mrk96} 
Markushevich, D., 
{\em Minimal discrepancy for a terminal cDV singularity is 
$1$}, {J. Math. Sci. Univ. Tokyo {\bf 3} (1996), no. 2, 445--456.}

\bibitem{MP04}
McKernan, J.; Prokhorov, Y.,
{\em Threefold thresholds.}  
{Manuscripta Math. 114 (2004), no. 3, 281--304.} 

\bibitem{MMM}
Mori, S.; Morrison, D. R.; Morrison, I.,
{\em On four-dimensional terminal quotient singularities.}
{Math. Comp. 51 (1998), no. 184, 769--786.}

\bibitem{Mircea02}
Musta\c{t}\v{a}, M.,
{\em Singularities of pairs via jet schemes.}
{J. Amer. Math. Soc. 15 (2002), no. 3, 599--615.}

\bibitem{Nadel91}
Nadel, A. M.,
{\em The boundedness of degree of Fano varieties 
with Picard number one.}  
{J. Amer. Math. Soc. 4 (1991), no. 4, 681--692.}

\bibitem{Nakayama04}
Nakayama, N.,
{\em Zariski-decomposition and abundance.}
{MSJ Memoirs, 14. Mathematical Society of Japan, 
Tokyo, 2004.} 

\bibitem{Oda88}
Oda, T.,
{\em Convex bodies and algebraic geometry. An introduction 
to the theory of toric varieties.} 
{Ergebnisse der Mathematik und ihrer Grenzgebiete (3) 15. 
Springer-Verlag, Berlin, 1988.}

\bibitem{Pro01}
Prokhorov, Y. G.,
{\em Lectures on complements on log surfaces.} 
{MSJ Memoirs, 10. Mathematical Society of Japan, Tokyo, 2001.}

\bibitem{Pukh04}
Pukhlikov, A. V.,
{\em Birational geometry of Fano direct products.}
{preprint math.AG/0405011.}

\bibitem{Reid79}
Reid, M.,
{\em Canonical $3$-folds.}  
{Journ\'ees de G\'eometrie Alg\'ebrique d'Angers, Juillet 
1979/Algebraic Geometry, Angers, 1979, pp. 273--310, 
Sijthoff \& Noordhoff, Alphen aan den Rijn---Germantown, 
Md., 1980.}

\bibitem{Sh85}
Shokurov, V. V.,
{\em A nonvanishing theorem.}
{Izv. Akad. Nauk SSSR Ser. Mat. 49 (1985), no. 3, 635--651.}

\bibitem{Sh88} 
Shokurov, V.V, 
{\em Problems about Fano varieties}, 
{Birational Geometry of Algebraic Varieties, 
Open Problems - Katata 1988, 30--32.}

\bibitem{Sh91} 
Shokurov, V.V, 
{\em A.c.c. in codimension $2$},
{preprint 1991.}

\bibitem{Sh92}
Shokurov, V. V.,
{\em Three-dimensional log perestroikas.
(Russian) With an appendix in English by Yujiro Kawamata.}  
{Izv. Ross. Akad. Nauk Ser. Mat. 56 (1992), no. 1, 105--203;  
translation in Russian Acad. Sci. Izv. Math. 40 (1993),  
no. 1, 95--202.}

\bibitem{Sh96}
Shokurov, V. V.,
{\em $3$-fold log models.}
{Algebraic geometry, 4., J. Math. Sci. 81 (1996), no. 3, 
2667--2699.} 

\bibitem{Sh00}
Shokurov, V. V.,
{\em Complements on surfaces. Algebraic geometry, 10.}  
{J. Math. Sci. (New York) 102 (2000), no. 2, 3876--3932.}

\bibitem{Sh04}
Shokurov, V. V.,
{\em Letters of a bi-rationalist. V. Minimal log discrepancies 
and termination of log flips.} 
{(Russian) Tr. Mat. Inst. Steklova 246 (2004), Algebr. Geom. 
Metody, Svyazi i Prilozh., 328--351;  
translation in  Proc. Steklov Inst. Math. 2004, no. 3 (246), 
315--336.} 

\bibitem{Siu95}
Siu, Y.-T.,
{\em The Fujita conjecture and the extension theorem 
of Ohsawa-Takegoshi.}  
{Geometric complex analysis (Hayama, 1995), 577--592, World 
Sci. Publ., River Edge, NJ, 1996.}

\bibitem{Viehweg95}
Viehweg, E.,
{\em Quasi-projective moduli for polarized manifolds.} 
{Ergebnisse der Mathematik und ihrer Grenzgebiete (3), 
30. Springer-Verlag, Berlin, 1995.}

\bibitem{Yasuda03}
Yasuda, T.,
{\em Dimensions of jet schemes of log singularities.}  
{Amer. J. Math. 125 (2003), no. 5, 1137--1145.}

\bibitem{Zariski62}
Zariski, O.,
{\em The theorem of Riemann-Roch for high multiples 
of an effective divisor on an algebraic surface.}  
{Ann. of Math. (2) 76 (1962), 560--615.}

\end{thebibliography}
\end{document}